\documentclass[11pt]{article}

\usepackage{a4wide}
\DeclareMathAlphabet\mathbb{U}{fplmbb}{m}{n}
\usepackage{color}
\usepackage{url}
\usepackage{amssymb}%
\usepackage{pstricks,pst-node,pst-text,pst-3d}%
\usepackage{graphicx}
\usepackage{amsthm}
\usepackage{psfrag}
\usepackage{floatflt}
\usepackage{algorithm}

\graphicspath{{Figures/}}
\input{psfig.sty}

\def\ni{\noindent}
\def\DD{\Delta}
\def\adj{\sim}

\def\adj{\sim}

\def\bs{\bigskip}

\def\cc{cc}
\def\Vt{V_{\geq 3}}
\def\nt{n_{\geq 3}}
\def\nongoob{n_G} 
\def\aa{\alpha}
\def\claim#1{\medbreak\ni{\bf Claim~#1}}
\def\term#1{{\em #1}\marginpar{\raggedright{\small\it #1}}}
\def\notat#1{{$#1$}\marginpar{\raggedright{\small $#1$}}}

\def\proof{\ni{\sl Proof.}\/\ }

\def\TEN{2-necklace}
\def\TEB{2-blossom}
\def\Gs{G_7}

\renewcommand{\qed}{\mbox{}\hspace*{\fill}{$\Box$}\medskip}
\def\qedclaim{\hfill$\triangle$\smallskip}
\def\LP{\mathcal{P}}
\def\int#1{I(#1)}
\def\ST{spanning tree}

\def\CVF{\overline{V(F)}}
\def\bs{\backslash}
\def\maxleaf{{\sc MaxLeaf}}
\def\H{S}
\def\q{\hspace{0.5cm}}

\newtheorem{theorem}{Theorem}
\newtheorem{repeatthm}{Theorem}
\newtheorem{proposition}{Proposition}
\newtheorem{lemma}{Lemma}
{\theoremstyle{definition} \newtheorem{definition}{Definition}}

\title{Spanning Trees with Many Leaves in Graphs without Diamonds and Blossoms}
\author{Paul Bonsma\thanks{Supported by the Graduate School ``Methods
    for Discrete Structures'' in Berlin, DFG grant GRK 
     1408} \hskip2cm Florian Zickfeld \thanks{Supported by the
     Studienstiftung des deutschen Volkes}\\[2mm]
\small Technische Universit\"at Berlin, Fachbereich Mathematik \\
\small Str. des 17. Juni 136, 10623 Berlin, Germany \\
\small \url{{bonsma,zickfeld}@math.tu-berlin.de}\\
}
\date{}

\begin{document}

\maketitle

\begin{abstract} 
It is known that graphs on $n$ vertices with minimum degree at least 3
have spanning trees with at least $n/4+2$ leaves and that this can be 
improved to $(n+4)/3$ for cubic graphs without the diamond $K_4-e$
as a subgraph.
We generalize the second result by proving that every graph with
minimum degree at least 3, without diamonds and certain
subgraphs  called blossoms, has a spanning tree with at least 
$(n+4)/3$ leaves, and generalize this further by allowing vertices of lower degree. 
We show that it is necessary to exclude blossoms in
order to obtain a bound of the form $n/3+c$. 

We use the new bound to obtain a simple FPT algorithm, which decides
in $O(m)+O^*(6.75^k)$ time whether a graph of size $m$ has a spanning tree
with at least $k$ leaves. This improves the best known time complexity
for {\sc Max Leaf Spanning Tree}.

\end{abstract}

\section{Introduction}

In this paper we study spanning trees with many leaves. We prove
a new extremal result, and apply it to obtain a fast
FPT algorithm for the related decision problem~\maxleaf.

We first introduce the extremal problem and explain our
contribution. Throughout this paper $G$ is assumed to be a
simple and connected graph on $n\geq 2$ vertices. Other graphs may be
multi-graphs, disconnected, or a $K_1$. The minimum vertex degree of
$G$ is denoted by~$\delta(G)$. Vertices of degree~1 are called leaves.

Linial and Sturtevant~\cite{LS87} and Kleitman and West~\cite{KW91}
showed that every graph $G$ with $\delta(G)\geq 3$
has a \ST\  with at least $n/4+2$ leaves, and that this bound is best
possible. The paper~\cite{KW91} also improves on this bound for graphs of
higher minimum degree. 

The examples showing that $n/4+2$ is best possible for graphs of
minimum degree~3 all consist of cubic  diamonds connected in a cyclic
manner. A {\em diamond} is the graph $K_4$ minus one edge, and an 
induced diamond subgraph of a graph $G$ is a {\em cubic diamond} if its
four vertices all have degree~3 in $G$, see
Figure~\ref{fig:diamond+blossom}~(a).

\PsFigCap{55}{diamond+blossom}{A cubic diamond (a), a \TEN~(b), and a
  \TEB~(c).}

Since these examples are very restricted it is natural to ask if
better bounds can be obtained when 
diamonds are forbidden as subgraphs. 
This question was answered by Griggs, Kleitman and  Shastri~\cite{GKS}
for {\em cubic} graphs, which are graphs where every vertex has
degree~3. They show that a cubic graph $G$ without diamonds
always admits a \ST\ with at least $n/3+4/3$ leaves.  For minimum
degree~3 the following bound is proved in~\cite{Bon06}. A graph $G$ with $\delta(G)\geq
3$, without cubic diamonds, contains a \ST\ with at least
$2n/7+12/7$ leaves. Both bounds are best possible for their 
respective classes.  

In~\cite{Bon06} it is conjectured that the following statement holds,
which would improve the bound $2n/7+12/7$ with only a minor extra
restriction, and would also generalize
the result for cubic graphs from~\cite{GKS}:
Every graph $G$ with $\delta(G)\geq 3$ and without {\em \TEN s}
contains a \ST\ with at least $n/3+4/3$ leaves. 
Informally speaking, a \TEN\ is a concatenation of $k\geq 1$
diamonds with only two outgoing edges, see Figure~\ref{fig:diamond+blossom}~(b).

In Section~\ref{sec:obstructions} we disprove this conjecture by constructing
graphs with $\delta(G)=3$ without \TEN s, which do not
admit \ST s with more than $4n/13+2$ leaves. On the positive side, we
prove that the statement is true after only excluding one 
more very specific structure, called a {\em \TEB}, see
Figure~\ref{fig:diamond+blossom}~(c).
Precise definitions of \TEN s and \TEB s are given in Section~\ref{sec:obstructions}.
So we prove that graphs $G$ with $\delta(G)\geq 3$ without \TEN s or \TEB s have a \ST\ with at least $n/3+4/3$ leaves. 
In fact we generalize this statement even further by removing any
restriction on the minimum degree. The resulting statement is given in Theorem~\ref{thm:main}, which is our main result.

Let
\notat{\Vt(G)} denote the set of vertices in $G$ with degree at
least~3 and \notat{\nt(G)} its cardinality. Let \notat{\ell(T)} be the
number of  leaves of a graph~$T$.  
\newcounter{thm:main}
\setcounter{thm:main}{\value{theorem}}
\begin{theorem}\label{thm:main}
Let $G$ be a simple, connected graph on at least two vertices 
which contains neither \TEN s nor \TEB s. Then, $G$ has a
spanning tree $T$ with
\[
\ell(T) \geq \nt(G)/3+
\left\{\begin{array}{ll}
4/3 & \mbox{ if } \delta(G)\geq 3 \\    
2   & \mbox{ if } \delta(G)\leq 2 \mbox{\/.} 
\end{array}\right.
\]
\end{theorem}
Section~\ref{sec:obstructions} shows that Theorem~\ref{thm:main} is also best possible for non-cubic graphs.
Without proof we remark that also other results can be generalized in
a similar way, e.g. it is not hard to extend the proof in~\cite{KW91}
to prove that all graphs $G$ have a \ST\ with at least
$(n-n_2(G))/4+c$ leaves, where $n_2(G)$ denotes the number of
vertices of degree~2 in~$G$.

Our proof of Theorem~\ref{thm:main} is constructive and can be turned
into a polynomial time algorithm for the construction of a spanning tree.
The main technical contribution of this paper is that we prove this
generalization of the statement in~\cite{GKS}, and improvement of the
statement in~\cite{Bon06}, without a proof as lengthy as the proofs in
these two papers. This is made possible by extending the techniques and
proofs from~\cite{GKS}. In Section~\ref{sec:mainproof} we argue that
the long case study in~\cite{GKS} actually proves a strong new lemma,
which we use as an important step in the proof of Theorem~\ref{thm:main}.  
We share the opinion expressed in~\cite{GKS} that a shorter
proof of the bound for cubic graphs might not exist. Therefore using
that result in order to prove the more general statement seems
appropriate.\medskip

We now explain the consequences that Theorem~\ref{thm:main} has for FPT
algorithms  (short for {\em fixed parameter
   tractable}) for the following decision problem.\medskip

\ni Max-Leaves Spanning Tree (\maxleaf):\\
INSTANCE: A graph $G$ and integer $k$.\\
QUESTION: Does $G$ have a spanning tree $T$ with $\ell(T)\geq k$? \medskip

It is known that \maxleaf~is $\mathcal{NP}$-complete, see~\cite{GJ79}.
When choosing $k$ as a parameter, an algorithm for this problem is
called an \term{FPT algorithm} if its complexity is bounded by
$f(k)g(n)$, where $g(n)$ is a polynomial. See~\cite{FG06} and~\cite{DF99} 
for introductions to FPT algorithms. $f(k)$ is called
the \term{parameter function} of the algorithm. Usually, $g(n)$ will
turn out to be a low degree polynomial, thus 
to assess the
speed of the algorithm it is mainly important to consider the growth
rate of $f(k)$. Though since \maxleaf\ is $\mathcal{NP}$-complete, $f(k)$ will most
likely always be exponential. Bodlaender~\cite{Bod93} constructed
the first FPT algorithm for \maxleaf\, with a parameter function of
roughly $(17k^4)!$. Since then, considerable effort has been put in finding
faster FPT algorithms for this problem, see
e.g.~\cite{DF95,FMR00,BBW03,EFL05,Bon06}.
The papers~\cite{BBW03,EFL05,Bon06} also establish a strong connection between extremal
graph-theoretic results and fast FPT algorithms. In~\cite{BBW03}, the
bound of $n/4+2$ from~\cite{KW91} mentioned above is used to find an
FPT algorithm with parameter function $O^*({4k\choose  k})\subset
O^*(9.49^k)$. Here the \notat{O^*} notation ignores polynomial
factors. With the same techniques the bound of $2n/7+12/7$ mentioned
above is is turned into the so far fastest algorithm, with a parameter
function in $O^*({3.5k\choose  k})\subset O^*(8.12^k)$
(\cite{Bon06}). Similarly  Theorem~\ref{thm:main} yields a new FPT
algorithm for \maxleaf, presented in Section~\ref{sec:FPT}. 
\newcounter{thm2}
\setcounter{thm2}{\value{theorem}}
\begin{theorem}\label{thm:FPT}
There exists an FPT algorithm for \maxleaf\ with time complexity
$O(m)+O^*(6.75^k)$, where $m$ denotes the size of the input graph and
$k$ the desired number of leaves.
\end{theorem}
This algorithm is the fastest FPT
algorithm for \maxleaf\ at the moment, both optimizing the dependency
on the input size and the parameter function. It simplifies the ideas
introduced by Bonsma, Brueggemann and Woeginger~\cite{BBW03} and is
also significantly simpler than the other 
recent fast FPT algorithms. Hardly any preprocessing of the input
graph is needed, since Theorem~\ref{thm:main} is already formulated
for a very broad graph class. We end in Section~\ref{sec:conclusions}
with a discussion of possible extensions  and further consequences of
Theorem~\ref{thm:main}. 

\section{Obstructions for Spanning Trees with Many Leaves}\label{sec:obstructions}

\subsection{Diamond Necklaces, Blossoms and Flowers}

As mentioned in the introduction \TEN s have been
identified as an obstruction for the existence of spanning trees with $n/3+c$
leaves in graphs with minimum degree~3, see~\cite{KW91}
and~\cite{Bon06}. In this section we show that 
they are not the only such obstruction. We start by precisely defining
\TEN s and \TEB s.  

The degree of a vertex $v$ in a graph $G$ is denoted by
$d_G(v)$ and by $d(v)$ if ambiguities can be excluded. 
A vertex $v$ of a subgraph $H$ of $G$  with $d_H(v)<d_G(v)$ is called
a \term{terminal} of $H$.

\begin{definition}[2-Necklace]\label{def:diamond} The graph
  $K_4$ minus one edge is called a \term{diamond} and denoted by
  $N_1$. The degree~3 vertices are the \term{inner vertices} of the 
  diamond.

  For $k\geq 2$ the \term{diamond necklace} $N_k$ is obtained from the
  graph $N_{k-1}$ and a vertex 
  disjoint $N_1$ by identifying a degree~2 vertex of $N_1$ with a
  degree~2 vertex of $N_{k-1}$. Thus, $N_k$ has again two degree~2 vertices,
  which are denoted by $c_1$ and $c_2$.
  
  An $N_k$ subgraph of $G$ is a \term{\TEN} if it only has $c_1$ and
  $c_2$ as terminals,   which both have degree~3 in $G$.
  See Figures~\ref{fig:diamond+blossom}~(a) and~(b). 
  If $G$ contains an $N_1$ this way, this $N_1$ subgraph is also
  called a \term{cubic diamond} of $G$. 
\end{definition}
Diamond necklaces will also be called {\em  necklaces} for short.
In the course of studying the leafy tree problem we found that the
subgraphs defined next are also an obstacle for the existence of
spanning trees with $n/3+c$ leaves.  

\begin{definition}[2-Blossom]
\label{def:blossom}
  The graph $B$ on seven vertices shown in
  Figure~\ref{fig:blossom+flowerX}~(a) is the blossom graph.
  A blossom subgraph $B$ of $G$ is a \TEB\ if $c_1$ and $c_2$ are its
  only terminals, and they both have degree~3 in $G$, see
  Figure~\ref{fig:blossom+flowerX}~(b). 
\end{definition}

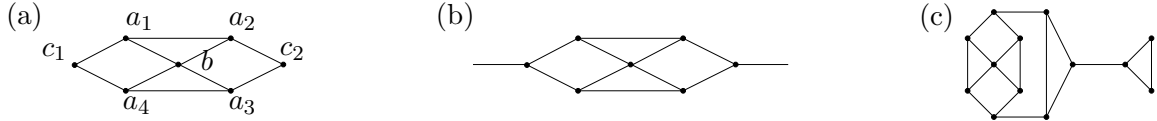
\PsFigCap{55}{blossom+flowerX}{A blossom graph~(a), a 2-blossom~(b), and a flower~(c).} 

If $G$ contains a 2-blossom $B$, only the
vertex $b$ has degree~4 in $G$, and the remaining vertices of $B$ have
degree~3 in $G$.  The two outgoing edges of a \TEN~respectively a
\TEB~may in fact be the same edge, in that case $G$ is just a
2-necklace respectively 2-blossom plus one additional edge.
The next proposition shows how many leaves can be gained within a blossom.
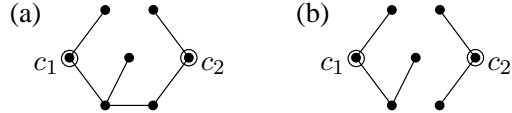
\PsFigCap{100}{BlossomTreeForms}{Spanning trees restricted to a blossom} 
\begin{proposition}\label{prop:blossomtrees}
Let $G$ be a graph with a blossom subgraph $B$ that has $c_1$ and $c_2$ as its only terminals. A \ST\ $T$ of $G$ exists with maximum number of leaves, such that $E(T)\cap E(B)$ has one of the forms in Figure~\ref{fig:BlossomTreeForms}.
\end{proposition}
\proof 
Consider a \ST\ $T$ of $G$ with maximum number of leaves. 
We may distinguish the following two cases for $E(T)\cap E(B)$: either this edge set induces a tree, or it induces a forest with two components, one containing $c_1$ and the other containing $c_2$. 

In the first case, at most three non-terminal vertices of $B$ can be
leaves of $T$, since a path from $c_1$ to $c_2$ contains at least two
internal vertices. In addition, if one of $c_1$ and $c_2$ is a leaf of
$T$, then $T$ can be seen to have at most two non-terminal vertices of
$B$ among its leaves. Since $c_1$ and $c_2$ together form a vertex cut
of $G$, one of them is not a leaf in $T$. It follows that replacing
$E(T) \cap E(B)$ by the edge set in
Figure~\ref{fig:BlossomTreeForms}~(a) does not decrease the number of
leaves. Since this edge set forms again a \ST\ of $B$, the resulting
graph is a \ST\ of $G$. 

Now suppose $E(T)\cap E(B)$ forms two components. At most four
non-terminal vertices of $B$ can be leaves of $T$. If one of $c_1$ and
$c_2$ is a leaf of $T$, then $T$ can have at most three non-terminal
vertices of $B$ among its leaves. One of $c_1$ and $c_2$ is not a leaf
in $T$, and thus it follows again that replacing $E(T) \cap E(B)$ by
the edge set in Figure~\ref{fig:BlossomTreeForms}~(b) does not
decrease the number of leaves, while maintaining a \ST\ of $G$.\qed 

We now present a family of graphs with minimum degree~3 which do not
contain diamond necklaces but do not have spanning trees with $n/3+c$
leaves. 

\begin{definition}[Flower, Flowerbed]\label{def:flower}   
  The \term{flower graph} is the graph on thirteen vertices shown in
  Figure~\ref{fig:blossom+flowerX}~(c). 

  The \term{flowerbed} $R_i$ of length $i$ consists of $i$ flowers,
  connected in a cyclic manner, see Figure~\ref{fig:flowerring2}. Formally,
  $R_i$ is constructed by starting with $i$ disjoint flowers, and
  adding $i$ edges in such a way that the graph is connected and has
  minimum degree 3. 
\end{definition}
Figure~\ref{fig:flowerring2} shows the flowerbed $R_5$. The solid edges
show a spanning tree with $4n/13+2$ leaves, which we will show to be optimal. 

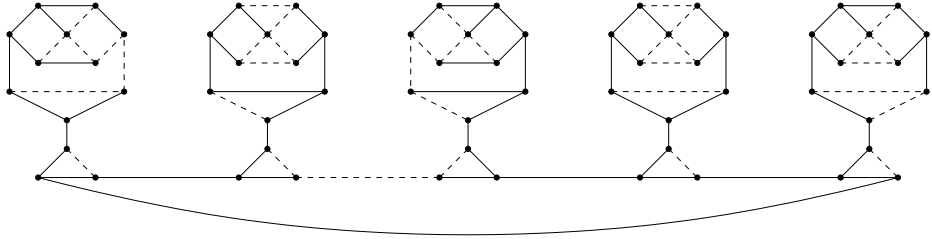
\PsFigCap{60}{flowerring2}{The flowerbed $R_5$. The solid edges show
  a tree with maximum number of leaves.}
\begin{proposition}\label{prop:flowerring2}
The flowerbed $R_i$ has no spanning tree with more than $4n/13+2$ leaves.
\end{proposition}
\proof 
Let $F$ be a flower in $R_i$ containing a blossom $B$, where $c_1$ and $c_2$ are the terminals of $B$. The neighbor of $c_j$ not in $B$ is called $f_j$ ($j=1,2$).
We will argue that no spanning tree $T$ of $R_i$ has more than four
leaves among $V(B)\cup
\{f_1,f_2\}$. Proposition~\ref{prop:blossomtrees} shows that without
loss of generality we may assume that $E(T)\cap E(B)$ has one of the
two forms in Figure~\ref{fig:BlossomTreeForms}. If it has the first
form, then one of $f_1$ and $f_2$ may be a leaf of $T$, but not both
since together they form a vertex cut of $R_i$. If $E(T)\cap E(B)$ has
the second form, then $f_1$ and $f_2$ are both cut vertices of $T$, 
so neither can be a leaf.

Now we consider the other vertices of $R_i$ that may be leaves in $T$.
Let $C$ be the cycle in $R_i$ that joins the $i$ flowers, that is, the
facial cycle of length $2i$ in Figure~\ref{fig:flowerring2}.  
Suppose $v\in V(C)$ is a leaf of $T$. In $G-v$, all vertices of $C$
except one are cut vertices, so $T$ may have at most one other vertex
of $C$ as a leaf. It follows that at most two vertices from $C$ can be
leaves in a spanning tree $T$ of $R_i$. The remaining vertices of
$R_i$ that we have not considered yet (two for every flower) are cut
vertices of $R_i$ and therefore not leaves of $T$. 

Summarizing, any \ST\ has at most two leaves in $C$, and at most four
additional leaves for every flower. The statement follows.\qed

\subsection{Tightness of the Bound}

The bound $n/3+4/3$ for cubic graphs is shown to be tight
in~\cite{GKS}. Infinitely many examples are given with no more than
$n/3+2$ leaves. On the other hand it is shown that there exists only
one graph that ensures that the additive term $4/3$ can not be
increased: the 3-dimensional cube $Q_3$, which has eight vertices and
only admits four leaves. 

Because the bound is best possible for cubic graphs, our bound is best
possible as well. But also graphs with arbitrarily many vertices of
higher and lower degree can be constructed which do not admit more
than $\nt/3+2$ leaves. Figure~\ref{fig:Tightness3}~(a) shows such an
example with many degree 2 and degree 4 vertices (which is closely
related to one of the examples from~\cite{GKS}). 

The reason that the additive term cannot be increased to 2 is again
only one example: Figure~\ref{fig:Tightness3}~(b) shows a graph on
$n=7$ vertices that only admits $4=\nt/3+5/3$ leaves. This graph will
be called $\Gs$ in the remainder. This graph is in fact a blossom
plus two edges; deleting any edge between two degree 4 vertices yields
a 2-blossom.An additive constant of $5/3$ is possible for non-cubic
graphs, but we will not prove this statement in this version of the
paper. 
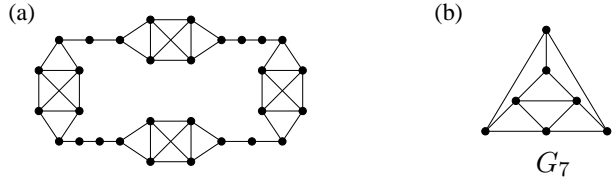
\PsFigCap{85}{Tightness3}{Non-cubic extremal graphs.}

\section{Proof of the Main Theorem}\label{sec:mainproof}

This section is devoted to the proof of the main theorem. We first
sketch the proof and give an overview of the different ingredients
that will be used. First we introduce a number of reduction rules in 
Section~\ref{subsec:reductions}. These reduction rules are applied to
the graph $G$ until an irreducible graph $G'$ is obtained. These rules
have the property that if the main theorem holds for every component
of $G'$, it also holds for $G$. In the next sections, we therefore
only have to consider irreducible graphs. In
Section~\ref{sec:usecubic} we argue that the proofs from~\cite{GKS}
for cubic graphs in fact show that if a non-spanning forest $F$ of
$G'$ is given, that contains all vertices of $G'$ of degree at least
4, then one of the trees of $F$ can be extended to a larger tree while
maintaining the proper leaf ratio. Finally, in
Section~\ref{sec:highdeg} we show how to obtain this starting forest
$F$ that covers all high degree vertices, while having enough
leaves. We use these tools in Section~\ref{subsec:wrapup} to prove
Theorem~\ref{thm:main}. 

\subsection{Reducible Structures}\label{subsec:reductions}

In this section we introduce a number of reduction rules. The proof of
the main theorem relies on locally extending a forest until it becomes
spanning while guaranteeing a certain number of leaves for every
intermediate forest. The reductions help to delay the treatment of
some substructures which cannot be readily handled during the
extension process and they also simplify the case study in the main proof. 

Ignoring rules that disconnect the graph,
the main idea behind the reduction rules is as
follows. A graph $G$ is reduced to a graph $G'$ with
$\nt(G)-\nt(G')=k$, such that every \ST\ of $G'$ can be turned into a
\ST\ of $G$  with at least $k/3$ additional leaves. This preserves the
desired leaf ratio. Lemma~\ref{reconstruction} states this idea more
precisely. 

We now give the necessary definitions. A vertex with degree at most two will be called a 
\term{goober}. We adopt this notion from~\cite{GKS}, although there it
is defined differently. In~\cite{GKS} goobers are those vertices of
degree at most two {\em resulting from a reduction rule}. We
observe that nowhere in the proofs the extra
structural information which this definition may provide is actually
used. Hence
goobers may simply be defined as we do here. The important gain is
that now we do not have to require graphs to have minimum degree~3 in our
statements. One important convention is that goobers are always
defined with respect to the whole graph, that is when we consider a subgraph 
$H$ of $G$, a vertex $v$ of $H$ is a goober if $d_G(v)\leq 2$. In our
figures, white vertices indicate goobers.   
A \term{high degree vertex} is a vertex of degree at least~4.

We first repeat the seven reduction rules defined in~\cite{GKS}, and
then introduce five new rules which 
are designed to handle structures containing higher degree
vertices. While the first seven rules are defined  in~\cite{GKS} for
graphs with maximum degree~3 we define them for arbitrary graphs, but
the vertices on which they act must have the same degrees as in the
original definition.

The seven reduction rules from~\cite{GKS} consist of graph operations
on certain structures, and conditions on when they may be applied. 
Figure~\ref{fig:GKSreductionsX} shows the operations. The black
vertices all have degree~3, and goobers are shown as white
vertices. Dashed edges are present in the resulting graph if and only
if they exist in the original graph.  The numbers above the arrows
indicate the decrease in $\nt$, and the numbers below the arrows
indicate the number of leaves that can be gained in a \ST\ when
reversing the reduction.

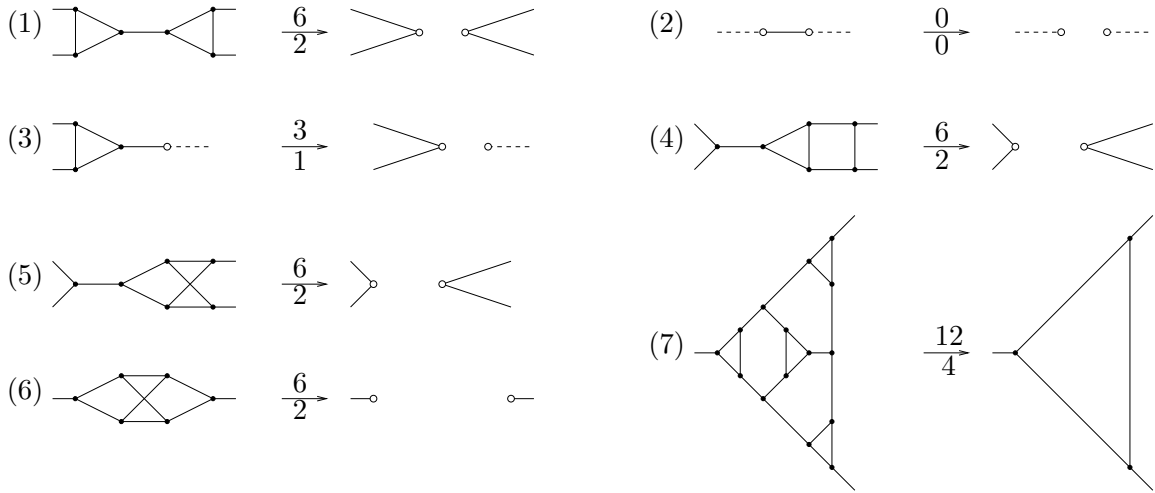
\PsFigCap{48}{GKSreductionsX}{The seven low degree reduction rules}  
\noindent The following restrictions are imposed on the application of
these rules (see Section~3 of~\cite{GKS}):  
\begin{itemize}
\item Reductions (1), (3), (4) and (5) may not be applied if the two
  outgoing edges from the left side, or the two outgoing edges from
  the right side, share a non-goober end vertex. (An outgoing edge
  from the left and an outgoing edge from the right may share a
  non-goober end vertex.)  
\item Reduction (7) may not be applied if any pair of outgoing edges
  shares an end vertex.  
\end{itemize}
In other words, a rule may not be applied if it would introduce
multi-edges incident with non-goobers, or if it would introduce a diamond.
These seven reduction rules will be called the \term{low-degree
  reduction rules}.  

We define an invariant that exhibits the properties
which should be maintained while doing graph reductions. 

\begin{definition}[Invariant]
A graph $H$ is said to \term{satisfy the invariant} if:
\begin{itemize}
\item $H$ is connected, or every component of $H$ contains a goober, and
\item every component of $H$ is either simple or it is a $K_2+e$, and 
\item $H$ contains neither \TEN s nor \TEB s.
\end{itemize}
\end{definition}
The reduction rules are applied in the induction step in the proof of
our main theorem; this invariant states the important properties that
should be preserved in the reduction process.

\newcounter{lem1}
\setcounter{lem1}{\value{lemma}}
\begin{lemma}
\label{cubic_reductions_invariant}
Let $G'$ be obtained from $G$ by the application of a low-degree
reduction rule. If $G$ satisfies the invariant then so does $G'$.
\end{lemma}
\proof 
Note that the reduction rules~(1)-(6) only introduce goobers as new
vertices and the only new edges are incident to these
goobers. Furthermore all other vertex degrees remain unchanged. Hence
these reductions cannot introduce \TEN s or \TEB s. Rule~(7) cannot
introduce a \TEB\  since a \TEB\ cannot 
share a vertex with a triangle induced by three vertices of degree
three. This is not true for \TEN s, but if rule~(7) introduces a \TEN,
two of the outgoing edges share an end vertex, contradicting
the condition for applying rule~(7). 

For all of the rules that may disconnect the graph, it is clear that
both new components will contain a goober. So the only way in which
one of the reductions might violate the invariant is by introducing
multiple edges. But using the imposed restrictions it can be seen that
multiple edges can only be introduced between two goobers, giving a $K_2+e$.\qed 

We now introduce five new reduction rules, which we call the
\term{high-degree reduction rules}. 
Each rule again consists of a graph operation and conditions on the
applicability. Figure~\ref{fig:reductions4} shows the graph
operations for the five rules.

The encircled vertices are the terminals,
which may have further incidences, unlike the other vertices. None of
the vertices in the figures 
may coincide, but there are no restrictions on outgoing edges sharing
end vertices. The numbers above the arrows indicate the decrease in $\nt$, and the
numbers below the arrows indicate the number of leaves that can be
gained in a \ST\ when reversing the
reduction.
Since (R4) must disconnect a component, this notion is not relevant
for (R4); this rule will be treated separately below. 
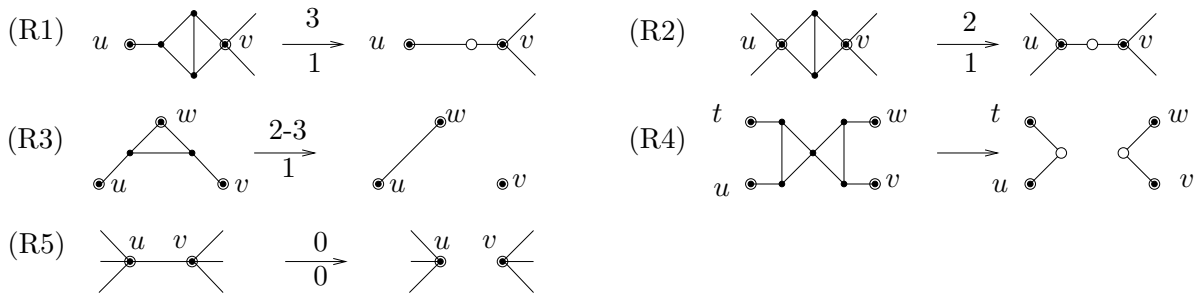
\PsFigCap{65}{reductions4}{The high-degree reduction rules.}

The following restrictions are imposed on the
applicability of these operations to a graph~$G$. First, none of the
reduction rules may be applied if it introduces a new \TEN or \TEB.
In addition, the following rule-specific restrictions are imposed. Let
\notat{cc(H)} denote the number of connected components of a graph
$H$.   
\begin{description}
\item[(R1)] $d_G(v)\geq 4$.
\item[(R2)] $d_G(u)\geq 4$ and $d_G(v)\geq 4$.
\item[(R3)] $cc(G')=cc(G)$, the edge $uw$ is not in $G$, and in addition
  $d_{G'}(v)\geq 3$, or $d_{G'}(w)\geq 3$, or both. 
\item[(R4)] $cc(G')>cc(G)$, that is $G'$ is {\em not} connected.
\item[(R5)] $d_G(u)\geq 4$, $d_G(v)\geq 4$, and $uv$ may not be a bridge.
\end{description}
A \term{bridge} is an edge whose deletion increases the number of
components. In the remainder, we will call a reduction rule
\term{admissible} if it can be applied without violating one of the
imposed conditions. 
In particular the condition that no \TEN s or \TEB s are introduced will be important.
Since the high-degree reductions rules are defined such that no \TEN s
or \TEB s can be introduced, it is easy to see that the following
lemma holds: 
\begin{lemma}
\label{highdeg_reductions_invariant}
Let $G'$ be obtained from $G$ by the application of a high-degree
reduction rule. If $G$ satisfies the invariant then so does $G'$.
\end{lemma}

Observe that in particular, (R5) seems counterproductive when the
goal is to find \ST s with many leaves, but it is useful to keep the
case analysis in the proof of Lemma~\ref{lem:start} simple. 

\begin{definition}[Reducible]
A graph $G$ is \term{reducible} if one of the low-degree
or high-degree reduction rules can be applied, and \term{irreducible} otherwise. 
\end{definition}

Griggs et al.~\cite{GKS} call a graph irreducible if none of the
low-degree reduction rules can be applied. Clearly, a graph that is
irreducible according to our definition is also irreducible according
to their definition, so we may apply their lemmas for irreducible
graphs also using the above definition of irreducibility.

Note that irreducible graphs satisfying the invariant are simple
because of reduction rule~(2). Components with only one vertex will be called
\term{trivial components} in the sequel. 

We now show that we can reverse all of these reduction rules while
maintaining spanning trees for every component, having the proper
number of leaves. For the low-degree reduction rules, this lemma was
implicitly proved in~\cite{GKS}. So for the detailed tree
reconstructions we refer to~\cite{GKS}, but we do repeat
the main idea behind the proof here.

\newcounter{lem2}
\setcounter{lem2}{\value{lemma}}
\begin{lemma}[Reconstruction Lemma] 
\label{reconstruction}
Let $G'$ be the result of applying a reduction rule to a connected
graph $G$ and $\aa\geq 0$. If $G'$ has $k$ non-trivial components $C_1,\ldots,C_k$,
which all have a spanning tree with at least $\nt(C_i)/3+\aa$ leaves, 
then $G$ has a spanning tree $T$ with
\[\ell(T)\geq\nt(G)/3+\aa k -2(k-1).\]
\end{lemma}

Note that the reduction rules create at most two components, that is
$k\leq2$. We use this lemma with $\aa=4/3$ if $G'$ is connected, and
with $\aa=2$  otherwise.\medskip

\proof
Suppose the applied rule was a low-degree reduction rule. Note that
$cc(G')$ is either 1 or 2. 
If $G'$ is connected, then its spanning tree can be turned into a
spanning tree of $G$ with $(\nt(G)-\nt(G'))/3$ more leaves. 
To prove this, it is shown in Section~3 of~\cite{GKS} for every rule
how to adapt the tree of $G'$ for $G$ ({\em tree
  reconstructions}). This already proves the statement if
$\cc(G')=1$. If $\cc(G')=2$, then applying the same tree
reconstructions yields a spanning forest of $G$ consisting of two
trees, with again $(\nt(G)-\nt(G'))/3$ more leaves in total. If both
components of $G'$ are non-trivial ($k=2$), then the two resulting
trees of $G$ can be connected to one \ST\ $T$ by adding one edge, losing
at most two leaves. In that case we have: 
\[
\ell(T)\geq \nt(G')/3+2\aa+(\nt(G)-\nt(G'))/3-2=\nt(G)/3+\aa k-2(k-1).
\]
If exactly one of the two components is trivial ($k=1$) then the
applied rule must be Rule (2) or (3). In this case, it can be checked
that after the tree reconstruction, one edge can be added without
decreasing the number of leaves; one leaf is lost but an isolated
vertex becomes a leaf. Then we have: 
\[
\ell(T)\geq \nt(G')/3+\aa+(\nt(G)-\nt(G'))/3=\nt(G)/3+\aa k-2(k-1).
\]
If both components of $G'$ are trivial ($k=0$), then the applied rule
was (2), and $G=K_2$, for which the statement holds: $-2(k-1)=2$. 
This proves the lemma when a low-degree reduction rule is applied. 

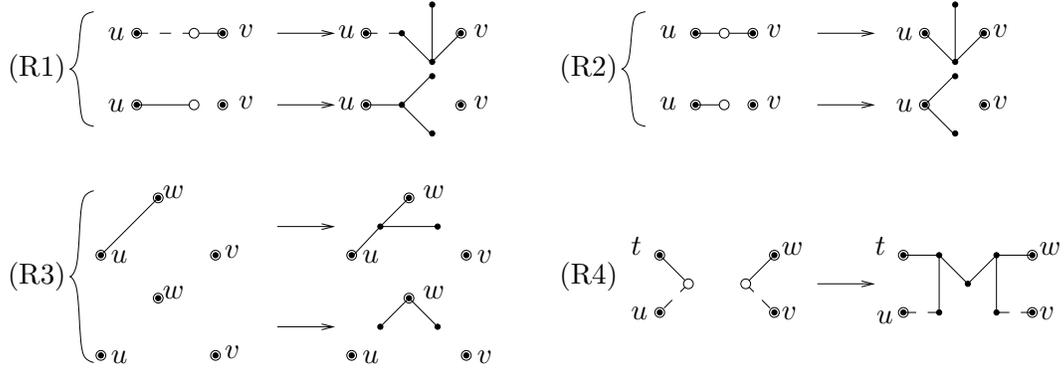
\PsFigCap{60}{reducreverseall}{Spanning tree constructions when
  reversing the new reduction rules} 

Now we consider the high-degree reduction rules. Note that rules (R1),
(R2), (R3) and (R5) do not increase the number of components, so $k=1$. So for
(R5) we do not have to change the \ST\ of $G'$. For (R1),
(R2) and (R3), Figure~\ref{fig:reducreverseall} shows how to gain at
least one additional leaf in every case, which suffices since each of
these rules decreases $\nt$ by at most three. Here it is 
essential that (R3) is admissible only if it creates at most one goober.
Dashed edges in the figure 
are present on the right if and only if they are present on the
left. Symmetric cases are omitted in the figure. Note that none of the
terminals of the operations can lose leaf status, except $w$ in the
second reconstruction for (R3). This is compensated by gaining two
new leaves here. So in every case enough leaves are gained to
maintain the ratio. 

Recall that (R4) is only admissible if it disconnects $G$ into two
components, which will be non-trivial, so
$k=2$. Figure~\ref{fig:reducreverseall} shows how 
to construct a \ST\ for $G$ from the two spanning trees for the
components, without decreasing the total number of leaves. Hence the
number of leaves of the resulting tree is at least 
\[
\nt(G')/3+2\aa=\nt(G)/3-5/3+2\aa>\nt(G)/3+2\aa-2=\nt(G)/3+\aa k-2(k-1)
\]
This proves the lemma for all reduction rules.\qed


The following property of irreducible graphs substantially simplifies
subsequent proofs. Here $\Gs$ denotes the graph from Figure~\ref{fig:Tightness3}~(b).

\newcounter{lem3}
\setcounter{lem3}{\value{lemma}}
\begin{lemma}[Edge Deletion]
\label{lem:highdeg}
Let $G$ be an irreducible graph not equal to $\Gs$
with adjacent vertices $u$ and $v$. If $d(u)=d(v)=4$, then $uv$ is a
bridge, or one of $u,v$ becomes an inner vertex of a cubic diamond
upon deletion of the edge $uv$. 
\end{lemma}
\proof
Suppose for the sake of contradiction a non-bridge edge $uv$
exists, between vertices of degree~4, such that none of $u,v$ becomes
an inner vertex of a diamond upon deletion of $uv$. 

Since $G$ is irreducible, no reduction rule is admissible. Clearly,
this must mean that a \TEN\ or \TEB\ is 
introduced when $uv$ is deleted, that is when (R5) is applied to $uv$. In
either case, we will derive a 
contradiction to the irreducibility of $G$. 

\claim{1} The graph $G-uv$ does not contain a \TEN\ $N$.\smallskip

\noindent Suppose for the sake of contradiction that $G-uv$ does contain a
\TEN\ $N$. Consider $N$ as a subgraph of $G$ (so $uv$ is counted
towards the degrees of $u$ and $v$). 

We first treat the case that $N$ consists of at
least two diamonds. If one of the diamonds in $N$ contains three
vertices of degree~3, we can use rule (R1), see
Figure~\ref{fig:edgedel1Y2}~(a). So now we may assume 
that one diamond on the end of the necklace contains $u$ as one of the
three vertices not shared with the next diamond, and the diamond on
the other end of the necklace contains $v$ this way. 

If $u$ is a vertex with degree~2 in $N$, then (R2) can be applied,
see Figure~\ref{fig:edgedel1Y2}~(b). This does not introduce a \TEN\
since the degree~4 vertex $v$ is part of $N$ on the other end. Because
$v$ is part of a diamond, this can also not introduce a \TEB. 
  
In the remaining case, both $u$ and $v$ are internal vertices of their
respective diamonds. Now it is admissible to apply (R5) to a different
edge incident with 
$u$, see Figure~\ref{fig:edgedel1Y2}~(c), where the dashed edge is the deleted one. 

This does not introduce a \TEN\ or \TEB: $u$ becomes part of a
triangle that is induced by degree~3 vertices, for which all outgoing
edges have different end 
vertices. Such a triangle cannot be part of a \TEB\ or \TEN. The other
end vertex of the deleted edge is still part of a diamond after
deletion, and thus is not part of a \TEB. It is not part of a \TEN\
since $v$ is in this part of the necklace.  

This concludes the case where $N$ consists of at least two diamonds.  
  
\PsFigCap{50}{edgedel1Y2}{Reductions when a long \TEN\ is created.}

Now suppose $N$ consists of a single diamond. If $u$ is an inner
vertex of this diamond, then $v$ cannot be part 
of the same diamond since we are dealing with simple graphs. This is then the
case we excluded by assumption, see Figure~\ref{fig:edgedel2Y}~(a). So
without loss of generality  $u$ is one of the vertices that have
degree~2 in the diamond. 

Now rule (R1) or (R2) is admissible, depending on whether $v$ is also
in the diamond, see Figures~\ref{fig:edgedel2Y}~(b) and~(c). 
This does not introduce a \TEB\ or \TEN, since in the case in
Figure~\ref{fig:edgedel2Y}~(b), a triangle containing a goober is
introduced, and in the case in Figure~\ref{fig:edgedel2Y}~(c), $v$ has
degree 4 and a goober at distance 2. 
Note that also no parallel edges are introduced: in the case in
Figure~\ref{fig:edgedel2Y}~(c) the edges leaving the diamond are distinct,
that is deleting $uv$ does not give a $K_4$, since in that case $uv$
would have been a bridge. This shows that it is admissible to apply
either (R1) or (R2), which contradicts the irreducibility of $G$.
\qedclaim  

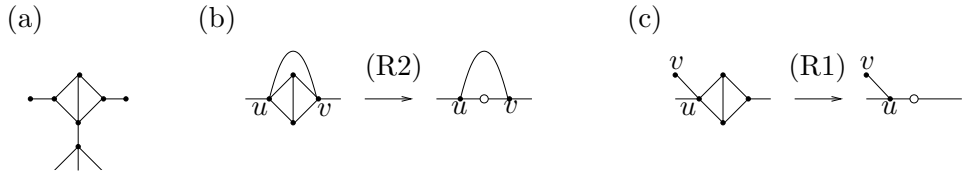
\PsFigCap{50}{edgedel2Y}{Reductions when a cubic diamond is created.}

\claim{2} The graph $G-uv$ does not contain a \TEB\ $B$.\smallskip

For $B$ we use the vertex labels from Figure~\ref{fig:blossom_ua2}~(a).
The degree~4 vertex of $B$ is labeled $b$, its terminals are called
$c$-vertices, and the remaining four vertices are called its
$a$-vertices. Now consider $B$ as a subgraph of $G$ (so $uv$ is
counted towards the vertex degrees). Since $d_G(u)=4=d_G(v)$ neither 
of them is equal to $b$, since $b$ has degree~4 even after the deletion of $uv$. 

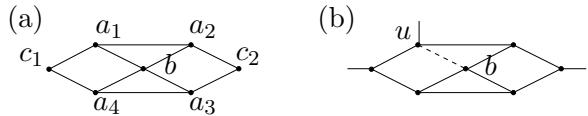
\PsFigCap{50}{blossom_ua2}{The blossom $B$ after deleting $uv$.}

If $u$ is an $a$-vertex, say without loss of generality $u=a_1$, then
it is admissible to delete the edge connecting $u$ to $b$ instead, see
Figure~\ref{fig:blossom_ua2}~(b). We argue that this does not
introduce a \TEB\ or \TEN. Figure~\ref{fig:blossom_edgedelcases} shows
the possible results of deleting $ub$ in more detail, depending on the
position of $v$.  
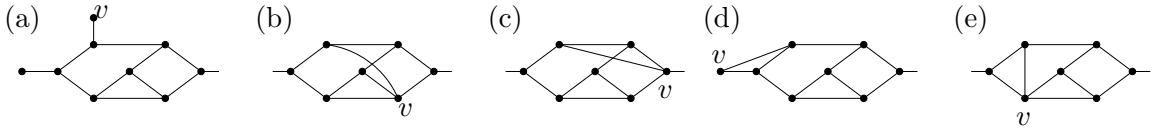
\PsFigCap{75}{blossom_edgedelcases}{Possible results of deleting $ub$.}
First suppose $v\not=a_4$.
After deleting $ub$, $b$ becomes part of a triangle that does not
share a vertex with another triangle, since we assumed $v\not=a_4$. It
follows that $b$ is neither part of a \TEN, nor of a \TEB. The vertex $u$
may be part of a triangle (when $v=c_2$ or when $v$ is not in $B$ but
adjacent to $c_1$), but such a triangle is not part of a diamond,
hence $u$ is not part of a \TEN. Finally we argue that $u$ is not part
of a \TEB: since $b$ is not part of a \TEB, its neighbor $a_2$ is not
part of a \TEB\ $B'$ unless it is a terminal of $B'$. In that case it
is not part of a triangle, but its neighbor $c_2$ is, which is
impossible. Hence $a_2$ is not part of a \TEB. Then if $u$ is part of
a \TEB\ $B'$, it must be a terminal of $B'$, and thus not part of a
triangle, but its neighbor $c_1$ must be part of a triangle. This is
again not possible. This concludes the proof that if $v\not=a_4$,
deleting $ub$ is an admissible application of (R5). 

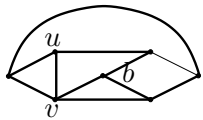
\PsFigCap{50}{blossomG_7X}{Deleting $ub$ yields a blossom.}
Now we need to consider the case case that $v=a_4$ and $u=a_1$, see
Figure~\ref{fig:blossom_edgedelcases}~(e). Deleting $ub$ does not
introduce a \TEN, but there is exactly one way in which it may
introduce a \TEB, which has $v$ as its central degree 4
vertex. Figure~\ref{fig:blossomG_7X} shows this case, the bold edges
indicate the new blossom. But now it can be seen that the original
graph, which includes $ub$, is exactly $\Gs$, a contradiction with our
assumption. We conclude that if $u$ is an $a$-vertex and $G\not=\Gs$,
in every case the edge $ub$ can be deleted by an admissible
application of (R5). 
  
It remains to consider the case that $u$ is a $c$-vertex. Then, (R3)
could be used, see Figure~\ref{fig:edgedel3X}~(a) and~(b). The bold
edges indicate the structure reduced by (R3). 
If in case (a) a \TEN\ is
introduced, $v$ would be an inner vertex of one of its diamonds, but
that is not possible since $d(v)=4$. In case (b) no \TEN\ can be
introduced, since $v$ is part of at most one triangle. 
In neither case a \TEB\ is introduced.\qedclaim 

\PsFigCap{35}{edgedel3X}{More reductions if a \TEB\ is created.}
We have thus derived a contradiction to the irreducibility of the
graph for all cases where deleting $uv$ would not be an admissible
application of (R5), which proves the lemma.\qed

\subsection{Using the Result for Cubic Graphs to Prove Tree Extendibility}\label{sec:usecubic}

Using our observation that the results in~\cite{GKS} hold when goobers
are simply defined as vertices with degree at most~2, we may restate
Theorem~3 from~\cite{GKS} as follows.  
%
\begin{theorem}\label{thm:GKS} Every irreducible graph $G$
  of maximum degree exactly 3 and without cubic diamonds
  has a spanning tree with at least $\nt(G)/3+\aa$
  leaves, where $\aa=4/3$ if $G$ is cubic and $\aa=2$ otherwise.
\end{theorem}
%

We give a short overview of the proof of this statement, as it appears
in~\cite{GKS}.  For a subgraph $T$ of $G$, in addition to
$\ell(T)$ the following values are considered.
By \notat{\nongoob(T)} we denote the number of non-goober vertices of
$G$ that are in $V(T)$. By \notat{\ell_d(T)} we denote the number of
\term{dead leaves} of $T$, that is leaves of $T$, which have no
neighbor in $V(G)\setminus V(T)$. 

The \term{value} of $T$ is defined as
$2.5\ell(T)+0.5\ell_d(T)-\nongoob(T)$. First it is shown in~\cite{GKS} that a tree $T$
with value at least 4 can always be found, and even one with value at
least 5.5 if $H$ contains at least one goober. Next, it is shown that
every non-spanning tree $T$ can be {\em extended}, that is a tree
supergraph $T'$ of $T$ can be found with value at least the value of $T$. 
This part of the proof consists of a rather involved case study. 
The extensions can be repeated until a spanning tree is found, in
which case all leaves are dead. Rewriting the value expression, and
rounding up the start value then yields Theorem~\ref{thm:GKS}. 

We observe that nowhere in the case study that proves extendibility
any information about the current tree $T$ is used; loosely speaking, 
only information about the part of $H$ `outside' of $T$ is used. In
particular, the fact that $T$ is connected is never used in the proof,
and neither are upper bounds on degrees of vertices already included
in $T$. 

We will now define the {\em leaf potential} of a subgraph which
generalizes the above definition of the value of a tree, and we will
formalize the notions `extendible' and `outside'. Using these new
notions a useful lemma can be formulated, which we conclude is proved,
but not stated in~\cite{GKS}. 

\begin{definition}[Leaf-Potential]
The \term{leaf-potential} of a subgraph $F\subseteq G$ is
\notat{\LP_G(F)}$=2.5\ell(F)+0.5\ell_d(F)-\nongoob(F)-6\cc(F).$
\end{definition}
If ambiguities are excluded in the context we simply write $\LP(F)$.

\begin{definition}[Extendible]
 Let $F$ be a subgraph of a graph $G$. Then $F$ is called \term{extendible} if there 
  exists an $F'$ with $F\subset F'\subseteq G$ and $\LP_G(F')\geq \LP_G(F)$.
\end{definition}

Above we already informally mentioned the subgraph of $G$ `outside' a
subgraph $F\subset G$. Considering the proof in~\cite{GKS}, we see that 
this graph may formally be defined as an edge induced graph as follows.
\begin{definition}[Graph Outside $\mathbf{F}$]
  Let $F$ be a non-spanning subgraph of $G$. 
  The \term{subgraph of $G$ outside of $F$} is
  \notat{F^C}$=G[\{uv\in E(G):u\not\in V(F)\}].$
  The \term{boundary} of $F$ is
  $
  V(F)\cap V(F^C)
  $
\end{definition}
Note that no edges between two vertices that are both in $V(F)$ are included in $F^C$.
If $G$ is clear from the context we call $F^C$ the \term{graph outside
  $F$}. Expressed using these definitions, the case study in~\cite{GKS}
yields the following lemma. 
\begin{lemma}[Extension Lemma]\label{lem:GKSextension}
Let $G$ be a connected
irreducible graph, and let $F\subset G$ such
that $F^C$ has maximum degree~3 and contains no cubic diamonds. Then
$F$ is extendible. 
\end{lemma}

\subsection{Growing Trees around High Degree Vertices}\label{sec:highdeg}

The purpose of this section is to prove Lemma~\ref{lem:start},
which grows trees around high degree vertices and yields a graph satisfying the
assumptions of Lemma~\ref{lem:GKSextension}.   Lemma~\ref{lem:start}
is the core of our proof of Theorem~\ref{thm:main}.

 We denote the set of vertices of $G$ which are not in $F$ by
$\CVF=V(G)\setminus V(F)$. The neighborhood $N(v)$ of a vertex $v$ is the set of all vertices
adjacent to $v$, and the closed neighborhood of $v$ is $N[v]=N(v)\cup 
\{v\}$. 
\term{Expanding} a vertex $v\in V(G)$, which is an operation on a subgraph $F$ of $G$,
yields a new subgraph with vertex set $V(F)\cup N[v]$, and edge set $E(F)\cup
\{uv:u\in N[v]\bs V(F)\}$. So all newly added neighbors of $v$ become
leaves, and $v$ may lose leaf status. The number of components
increases by one if and only if $v\not\in V(F)$. 
Expanding a list of vertices means expanding the vertices in the given order.

We adopt the short-hand notation
\notat{\DD(x,y,z)}$:=2.5y+0.5z-x$ from~\cite{GKS} to express the change in 
$\LP_G$ when extending a graph $F$ to a new graph $F'$.  
Let $\DD\ell$ denote $\ell(F')-\ell(F)$, and define $\DD\ell_d$ and
$\DD \nongoob$ analogously. So when $F$ and $F'$ have the same number of
components, the extension is valid if and only if 
$\DD(\DD\nongoob,\DD\ell,\DD\ell_d)\geq 0$, and if a new component is
introduced we need $\DD(\DD\nongoob,\DD\ell,\DD\ell_d)\geq6$. 

For the sake of simpler notation, instead of writing e.g. $\DD(\DD
\nongoob,\DD \ell,\DD \ell_d)\geq \DD(4,3,1)=4$, we will simply write
$\DD(4,3,1)=4$. Hence the three parameter values need not be 
exactly $\DD \nongoob$, $\DD \ell$ and $\DD\ell_d$ but reflect the
worst case scenario. That is, the change in the leaf potential that we
prove is to be read as a lower bound for the actual change.

\newcounter{lem5}
\setcounter{lem5}{\value{lemma}}
\begin{lemma}[Start Lemma]\label{lem:start} Let $G$ be an 
irreducible graph not equal to $\Gs$ and $F$ a (possibly empty) subgraph of $G$, such
that $F^C$ contains at least one vertex of degree at
least 4, and contains neither \TEN s or \TEB s. Then, $F$ is
extendible.
\end{lemma}

\proof
First suppose $F$ is not the empty graph. If there is a vertex $v$ on
the boundary of $F$ which is not a leaf,  then $F'$ can be obtained by
expanding $v$. There is no leaf lost since $v$ was not a leaf, and the
newly added vertices are leaves. So the augmentation inequality is
satisfied: $\DD(k,k,0)\geq0$. Hence, we may assume in the remainder
that only leaves of $F$ have neighbors in $\CVF$, or in other words,
all vertices on the boundary of $F$ are leaves of $F$. 

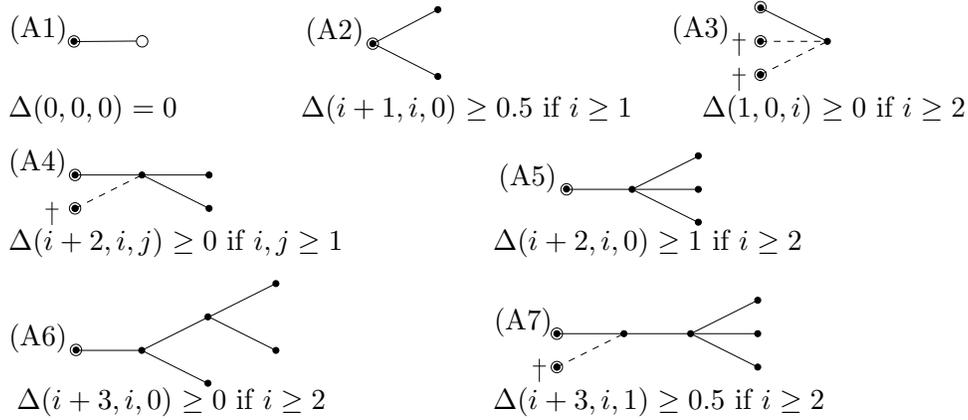
\PsFigCap{70}{AugmentationRules_vertcompr}{Simple augmentations of an existing subgraph}

The next step is the attempt to augment $F$ using the operations
(A1)-(A7), see Figure~\ref{fig:AugmentationRules_vertcompr}. Conventions for
this figure are that encircled vertices belong to $V(F)$, solid edges
show the expansion and vertex degrees shown are to be  understood
as lower bounds. Dead leaves are marked with a cross. All of the
expansions in the figure extend $F$ without creating a new
connected component, and satisfy
$\DD(\DD\nongoob,\DD\ell,\DD\ell_d)\geq0$. Thus the resulting graph
$F'$ is an extension as claimed in the lemma. Together these
augmentation rules yield the following claim.   

\claim{0} The subgraph $F$ is extendible, if a vertex in $V(F)$ has a
goober neighbor in $\CVF$ or at least two neighbors in $\CVF$, or if there 
is a high-degree vertex $v\in \CVF$ at distance at most two from~$F$.\medskip  

\ni
If a goober from $\CVF$ is adjacent to $F$ then (A1) can be applied.
If a vertex in $V(F)$ has at least two neighbors in $\CVF$, (A2) can
be applied. So from now on we will assume every vertex in $V(F)$ has
at most one neighbor in $\CVF$, and this neighbor is not a goober. If
a high-degree vertex in $\CVF$ is 
adjacent to a vertex in $V(F)$, (A3), (A4) or (A5) can be applied. The
creation of the dead leaves in (A3) and (A4) follows from the fact
that (A2) cannot be applied anymore. If a high-degree vertex in $\CVF$ has
distance two from a vertex in $V(F)$, (A6) or (A7) can be applied. \qedclaim 

The rest of the proof will handle the more complicated case when $F$ is the empty graph, or the
only high-degree vertices in $\CVF$ are at a larger distance from
$V(F)$. We then introduce a new component for $F$. This is more
complicated because adding a further component comes at a certain 
cost, more precisely we need that the new component satisfies
$\DD(\DD\nongoob,\DD\ell,\DD\ell_d)\geq6$.   

The rest of the proof is divided into three more claims. The
first one handles the easiest cases, and the second one handles all
cases except those where every degree~4 vertex is the common vertex
of two edge-disjoint triangles. This final case is then taken care of
in the third claim.  
Throughout the proof we assume, sometimes implicitly, that none of the
situations that have been handled earlier can occur.  

\claim{1} Let $v\in \CVF$, $d(v)\geq 4$, and $w\in N(v)$. In the
following four situations $F$ is extendible: $d(v)\geq5$, or
$d(v)=4$ and $w$ is a goober, or $d(v)=d(w)=4$, or $d(v)=4$ and $N[w]\subset 
  N[v]$.\medskip

\ni First note that no vertex in $N[v]$ or $N[w]$ is part of $F$ by Claim~0.
If $d(v)\geq 5$, expanding $v$ yields $\DD(k+1,k,0)\geq 6.5$, since
$k\geq 5$. For $d(v)=4$ and $w$ a goober, expanding $v$ gives $\DD(4,4,0)=6$. 
 
Now suppose $d(v)=d(w)=4$. If $vw$ is a bridge, expanding $v$ and $w$
yields $\DD(8,6,0)=7$. Otherwise, Lemma~\ref{lem:highdeg} shows that
either $v$ or $w$, say $v$, becomes the inner 
vertex of a cubic diamond upon deletion of the edge $vw$. (Note that
we assumed $G\not=\Gs$, so Lemma~\ref{lem:highdeg} may be applied.) Thus, either $w$ has
two neighbors not in $N[v]$ and expanding $v,w$ yields
$\DD(7,5,1)=6$ (see Figure~\ref{fig:StartLemClaim1X}~(a)), or $w$
shares two neighbors with $v$ in which case expanding $v$ yields
$\DD(5,4,3)=6.5$ (see Figure~\ref{fig:StartLemClaim1X}~(b)).

So now we may assume that all neighbors of $v$ have degree~3.
If $N[w]\subset N[v]$ then either the unique vertex $u\in N[v]\bs N[w]$ has
two neighbors not in $N[v]$, in which case expanding $u,v$ gives
$\DD(7,5,1)=6$, see Figure~\ref{fig:StartLemClaim1X}~(c), or there is
another vertex $x\in N[v]-w$ with $N[x]\subset N[v]$, and $v$ is
expanded to obtain $\DD(5,4,2)=6$, see Figure~\ref{fig:StartLemClaim1X}~(d).\qedclaim

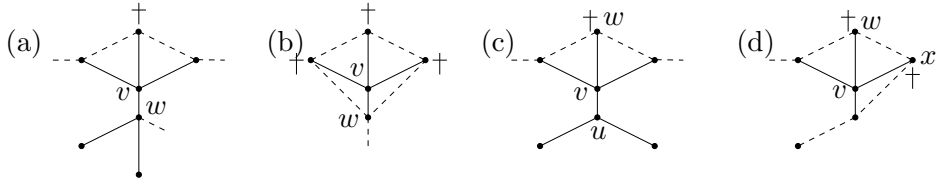
\PsFigCap{60}{StartLemClaim1X}{Figures for Claim 1.}

\noindent Summarizing, we may now assume that $\CVF$ contains no
vertices of degree at least 5, and if it contains a vertex $v$ of
degree~4, all neighbors of $v$ have degree~3 and have either one or two
neighbors not in~$N[v]$.\bigskip

\claim{2} If $\CVF$ contains a
  vertex $v$ with $d(v)=4$ and a vertex $w\in N(v)$ which has two neighbors
  $a,b\not\in N[v]$, then $F$ is extendible. \medskip

\ni We denote the other three neighbors of $v$ by $x,y,z$.  If one of
$a,b,x,y,z$ has all of its neighbors in $\{a,b\}\cup N[v]$, we have
$\DD(7,5,1)=6$ by expanding $v,w$, see
Figure~\ref{fig:StartLemClaim2aX}~(a). 

If $a$ or $b$ is a goober we obtain
$\DD(6,5,0)\geq 6.5$, see Figure~\ref{fig:StartLemClaim2aX}~(b). If $a$
or $b$ is adjacent to a vertex $c\in V(F)$, then expanding $v,w$ will
make $c$ a dead leaf and yields $\DD(7,5,1)\geq 6$, see
Figure~\ref{fig:StartLemClaim2aX}~(c). 
If one of $a,b,x,y,z$ has at least two neighbors not in
$N[v]\cup\{a,b\}$, we obtain $\DD(9,6,0)=6$ by expanding $v,w$ and
this vertex, see Figure~\ref{fig:StartLemClaim2aX}~(d).\medskip  

\PsFigCap{80}{StartLemClaim2aX}{Figures for Claim 2.}

Hence we may assume that $a,b,x,y,z$ each have exactly one
neighbor outside $N[v]\cup\{a,b\}$. This neighbor is not 
part of $F$. Since they all have degree at least 3, these
five vertices must induce three edges. This 
implies that one of $a,b$ has degree~4 since we already know that
$x,y,z$ have degree~3. We may assume without loss of generality
that $d(a)=4$ and $d(b)=3$. We  
distinguish two cases depending on whether $a$ is adjacent to $b$ or
not. The statement `$a$ is adjacent to $b$' is denoted by \notat{a\adj
  b}. We denote the neighbor of $x$ outside of $N[v]\cup\{a,b\}$ by
$x'$, and similarly  $a',b',y',z'$ are defined.\medskip    

\noindent\textbf{Case 1.} $a$ is adjacent to $x$ and $y$ while $b$
is adjacent to $z$.\smallskip 

Consider expanding $v,x,z$. All vertices in $\{a,b,w,y,x',z'\}$ are
adjacent to at least one of $v,x,z$, thus we have $\DD(9,6,1)=6.5$
unless $x'=z'$, see 
Figure~\ref{fig:StartLemClaim2bX}~(a). By an analogous argument with $y$ in
the place of $x$ we may now assume that $x'=z'=y'$. Then, expanding $v,x$
yields $\DD(7,5,1)=6$, since $y$ becomes a dead leaf, see
Figure~\ref{fig:StartLemClaim2bX}~(b).\medskip   

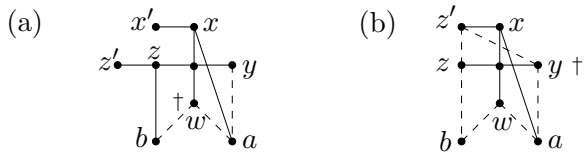
\PsFigCap{80}{StartLemClaim2bX}{Figures for Claim 2, Case 1}

\noindent\textbf{Case 2.}  $a$ is adjacent to $b$ and $x$ while $y$ is
adjacent to $z$.\smallskip 

If $x'\neq a'$, expanding $a,x,v$ yields $\DD(9,6,1)=6.5$, see
Figure~\ref{fig:StartLemClaim2cX2}~(a), so we may  assume that
$x'=a'=:c$, 
and this creates a situation symmetric in $b$ and $c$. 
By Claim~1 we have that $b\not\adj c$. Now first suppose $b'\adj c$. 
Then expanding $b',c,x,v$ yields $\DD(10,6,3)=6.5$, provided $b'$ has
a neighbor $d$ other than $y,z$, see
Figure~\ref{fig:StartLemClaim2cX2}~(b). Note that $d\in V(F)$ is not
possible since augmentation (A6) could have been applied instead.  
 
\PsFigCap{60}{StartLemClaim2cX2}{Figures for Claim 2, Case 2}

If $N(b')=\{b,c,y\}$, then (R3) is admissible, see
Figure~\ref{fig:StartLemClaim2cX2}~(c). Since $b'$ becomes a goober this
cannot introduce a \TEN. 
Hence it must be that $b'\adj y,z$ and the graph has $\DD(9,5,5)=6$, see
Figure~\ref{fig:StartLemClaim2dY}~(a). This concludes the cases with $b'\adj c$.

The case $b'\not\adj c$ can be excluded because then (R3) would be
admissible, see  Figure~\ref{fig:StartLemClaim2dY}~(b). Note that this
cannot create a \TEN\ involving $y,z$ since then (R2) would have been
admissible. 

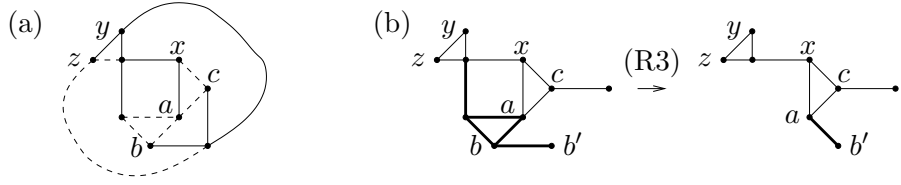
\PsFigCap{60}{StartLemClaim2dY}{Additional figures for Claim 2, Case 2.}

This concludes the proof of Claim~2.\qedclaim\medskip

Summarizing Claims~0,~1, and~2, we may now assume that all neighbors of a
degree~4 vertex $v\in \CVF$ have degree~3, and have exactly one
neighbor not in $N[v]$. In other words, $v$ is the common vertex of
two edge-disjoint triangles, see
Figure~\ref{fig:bowtie}. 

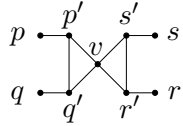
\PsFigCap{60}{bowtie}{The bow tie subgraph}

\claim{3} If the graph outside $F$ contains a vertex $v$ with
  $d(v)=4$ such that all its neighbors have degree~3 and one neighbor outside
  $N[v]$, then $F$ is extendible. \medskip

We denote the neighbors of $v$ by $p',q',r',s'$ and assume that $p'\adj
q'$ and $r'\adj s'$. The neighbor of $p'$ outside $N[v]$ is denoted by
$p$ and similarly $q,r,s$ are defined, see
Figure~\ref{fig:bowtie}. We split the proof of the claim 
into three cases. \medskip

\noindent \textbf{Case 1.} $p=q$\smallskip

If $p$ has degree~3, we can apply (R1), see
Figure~\ref{fig:StartLemClaim3aX}~(a). So now without loss of
generality $p$ has degree 
4. Then by Claim~2, $p$ is also part of two edge-disjoint
triangles. So if $p=r$ then also $p=s$. In that case we can expand
$p', v$ to obtain $\DD(6,4,4)=6$, see
Figure~\ref{fig:StartLemClaim3aX}~(b). So now $p\ne r$, $p\ne
s$. Consider applying (R2) to the diamond consisting of $p,p',q',v$, see
Figure~\ref{fig:StartLemClaim3aX}~(c). If this introduces a \TEN,
(R1) could have been applied to the diamond on the other end
of this necklace. It cannot introduce a \TEB\ since the triangles of
a \TEB\ contain a degree~4 vertex. \medskip

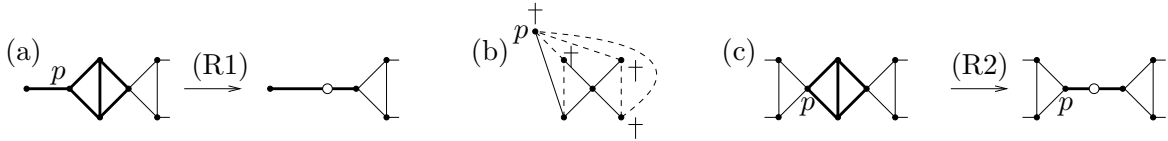
\PsFigCap{60}{StartLemClaim3aX}{Figures for Claim 3, Case 1}

\noindent \textbf{Case 2.} $p=r$\smallskip

Note that $d(p)=3$ by Claim~2. Also $q\neq s$ since the graph does not
contain \TEB s and the case that $d(q=s)=4$ is again excluded by
Claim~2. Thus, (R3) is admissible, see
Figure~\ref{fig:StartLemClaim3bX}. 
\medskip

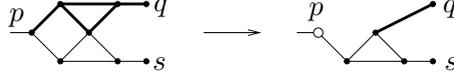
\PsFigCap{60}{StartLemClaim3bX}{The figure for Claim 3, Case 2}

\noindent\textbf{Case 3.} $p,q,r,s$ pairwise different.\medskip

In this case either (R4) or (R3) is admissible: if $v$ is a cut
vertex, (R4) may be used (it increases the number of components).  
Otherwise, we may assume  without loss of generality that $p$ and $s$
are in same connected component of $G-v$, and (R3) can be applied
without disconnecting the graph. See Figures~\ref{fig:StartLemClaim3cX}~(a) and~(b).\qedclaim

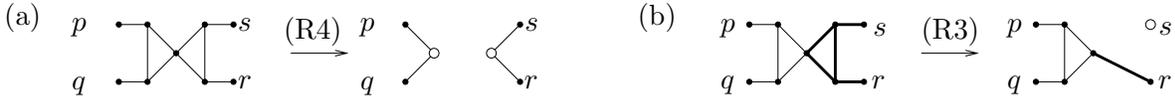
\PsFigCap{60}{StartLemClaim3cX}{Figures for Claim 3, Case 3}

This concludes all possible cases: whenever the subgraph of $G$
outside of $F$ contains a high degree vertex, we have shown that $G$
is either reducible, or $F$ is extendible.\qed

\subsection{The Proof of the Main Result}
\label{subsec:wrapup}

This section is devoted to combining the tools developed in the last
three subsection in order to prove Theorem~\ref{thm:main}, which we
repeat here for convenience.
\setcounter{repeatthm}{\value{thm:main}}
\begin{repeatthm}
Let $G$ be a simple, connected graph on at least two vertices 
which contains neither \TEN s nor \TEB s. Then, $G$ has a
spanning tree $T$ with
\[
\ell(T) \geq \nt(G)/3+
\left\{\begin{array}{ll}
4/3 & \mbox{ if } \delta(G)\geq 3 \\    
2   & \mbox{ if } \delta(G)\leq 2 \mbox{\/.} 
\end{array}\right.
\]
\end{repeatthm}
\proof We prove the statement by induction. For our induction
hypothesis we actually prove that the above 
statement holds for every connected graph which satisfies the
invariant. Then the statement follows for simple graphs.  
 
First suppose $G$ is irreducible. If $G$ has maximum degree exactly~3,
Theorem~\ref{thm:main} follows immediately from Theorem~\ref{thm:GKS}.
If $G$ has maximum degree at most~2, $G$ has a \ST\ with at least two leaves (note
that we assumed that $G$ is not a $K_1$), which suffices. 
If $G=\Gs$, then a \ST\ with $4=\nt(G)/3+5/3$ leaves can be obtained.
So we may now assume that $G$ contains at least one high degree vertex, and is not equal to $\Gs$. 

We start with an empty subgraph $F$ of $G$, which has $\LP_G(F)=0$. The
Start Lemma (Lemma~\ref{lem:start}) shows that, as long as there is at
least one high degree vertex not in $F$, we can extend $F$ while
maintaining $\LP_G(F)\geq 0$. When all high degree vertices are included
in $F$, the Extension Lemma (Lemma~\ref{lem:GKSextension}) can be
applied iteratively, until a spanning subgraph $F'$ is obtained with
$\LP_G(F')\geq 0$. Without loss of generality, we may assume that $F'$
is a forest; cycles can be broken without decreasing the number of
leaves. Since all leaves of a spanning subgraph are dead
we deduce
\[0\leq \LP_G(F')=3\ell(F')-\nt(G)-6\cc(F')\;
\Longrightarrow\;\ell(F')\geq \nt(G)/3+2\cc(F').\] 
We can now add $\cc(F')-1$ edges to
$F'$ to obtain a \ST, losing at 
most $2(\cc(F')-1)$ leaves, so the resulting tree has at least
$\nt(G)/3+2$ leaves.%

It remains to consider the case that $G$ is reducible (the induction
step). Some reduction rule is admissible, and the reduced graph
$G'$ again satisfies the invariant, by
Lemmas~\ref{cubic_reductions_invariant} and~\ref{highdeg_reductions_invariant}.

First suppose $G'$ is connected. None of the reduction rules remove
goobers, so if $\delta(G)\leq 2$, then $\delta(G')\leq 2$, and by
induction $G'$ has 
a \ST\ with at least $\nt(G')/3+2$ leaves. Lemma~\ref{reconstruction}
then shows that  $G$ has a \ST\ with at least $\nt(G)/3+2$
leaves. Similarly, if $\delta(G)\geq 3$ then it follows that $G$ 
has a \ST\ with at least $\nt(G)/3+4/3$ leaves.
Now suppose the reduction rule yields a disconnected graph $G'$. Then,
by the definition of the reduction rules, every resulting
component has a goober. So by induction, every non-trivial
component $C$ of $G'$ has a \ST\ with at least $\nt(C)/3+2$
leaves. Thus Lemma~\ref{reconstruction} implies that $G$ has a \ST\ with
at least $\nt(G)/3+2$ leaves.\qed

\section{A fast FPT Algorithm for \maxleaf}
\label{sec:FPT}

In this section we present a fast and relatively simple FPT algorithm
for \maxleaf, which uses Theorem~\ref{thm:main} as an essential
ingredient. The other two ingredients are a short preprocessing step,
consisting of two reduction rules, and an enumerative procedure, which
is similar to the one introduced in~\cite{BBW03}, and also
applied in~\cite{Bon06}. 

We start by presenting the two reduction rules that constitute the
preprocessing phase. The reduction rules for the FPT algorithm are
different from the rules used in Section~\ref{sec:mainproof}. It is
important that they yield an equivalent instance of the decision
problem \maxleaf, but there are no conditions on the ratio between
the decrease in vertices and in possible leaves. 

Recall that in  \TEN s and \TEB s both
terminals have degree~3 in $G$. The rules we introduce now
also reduce diamonds and blossoms whose two terminals have arbitrary
degree. However the two terminals of the subgraph must still be 
the two vertices that have degree~2 in the diamond
necklace or blossom itself. Such a subgraph of $G$ will be called a
\term{2-terminal diamond} respectively a \term{2-terminal blossom}.
Rule (F1) in Figure~\ref{fig:FPTreductionsX}, which resembles rule
(R2), reduces 2-terminal diamonds.  Since the \TEN\ $N_k$ consists
of $k$ 2-terminal diamonds it is reduced as well by rule~(F1).
Rule (F2) in Figure~\ref{fig:FPTreductionsX} reduces
2-terminal blossoms. The next lemma proves the correctness of these
rules.

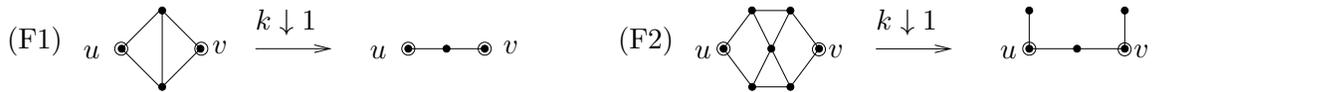
\PsFigCap{80}{FPTreductionsX}{Two reduction rules for an instance $(G,k)$ of \maxleaf.}
\newcounter{lem6}
\setcounter{lem6}{\value{lemma}}
\begin{lemma}
\label{FPT_reds_correct}
Let $G'$ be the result of applying reduction (F1) or (F2) to
$G$. Then $(G',k-1)$ is a YES-instance for \maxleaf\ if and only if $(G,k)$
is a YES-instance for \maxleaf. 
\end{lemma}
\proof
 First consider the case where the applied rule was (F1), reducing a
 2-terminal diamond $D$ of $G$ with terminals $u$ and $v$. 
 Consider a \ST\ $T$ of $G$ with at least $k$ leaves. 
 Observe that in $T$, we can always replace the set of edges $E(T)\cap
 E(D)$ by one of the two sets shown on the right in
 Figure~\ref{fig:FPTred_reconstr}~(a), or a symmetric set, 
 without decreasing the number of leaves. Then by replacing it by the
 corresponding structure on the left, we obtain a \ST\ $T'$ of $G'$
 with at least $k-1$ leaves. Note here that the terminals $u$ and $v$
 remain leaves in $T'$ if they are leaves in $T$. Similarly, any \ST\
 of $G'$ has one of the forms shown on the left when restricted to the
 two edges resulting from the reduction. Replacing it with the
 corresponding structure on the right shows that if $(G',k-1)$ is a
 YES-instance, $(G,k)$ is as
 well. Figure~\ref{fig:FPTred_reconstr}~(b) can be used to prove the
 statement for rule (F2) analogously. For this it is useful to note
 that Proposition~\ref{prop:blossomtrees} shows that we may assume
 that $T$ restricted to the blossom has one of the forms on the right
 of Figure~\ref{fig:FPTred_reconstr}~(b).\qed   

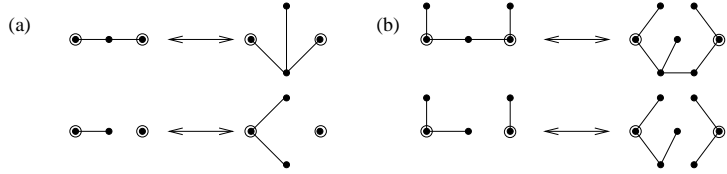
\PsFigCap{70}{FPTred_reconstr}{Tree reconstructions for (F1) and (F2).}

Throughout this section we will denote the set of leaves of a graph
$G$ by \notat{L(G)}. 
We now explain how to obtain a graph \notat{\H(G)} from a 
graph $G$ by suppressing vertices. \term{Suppressing} a
vertex $u$ of degree~2 means deleting $u$ and adding an edge between the two 
neighbors of $u$, or a loop if both edges incident with $u$ end in the same vertex.
We allow this operation to
introduce parallel edges and loops, so the degrees of non-suppressed
vertices are maintained. 
If $|\Vt(G)|=0$, that is $G$ is a path or cycle, then $\H(G)$ is the
empty graph.  If $|\Vt(G)|>0$ then $\H(G)$ is obtained from $G$ by 
suppressing all degree~2 vertices. 
So $V(S)=L(G)\cup \Vt(G)$, and $G$ is a subdivision of $\H(G)$. 
Hence loops and non-loop edges of $\H(G)$ correspond to cycles and
paths of $G$ respectively. 
Let $uv$ be a non-loop edge of $\H(G)$ where the corresponding path
$P_{uv}$ in $G$ has $i$ internal vertices. We define a cost function
$c$ on the non-loop edges of $\H(G)$ which assigns cost
\notat{c(uv)}$=\min\{i,2\}$ to $uv$.  
Thus $c(uv)$ is the maximum possible number of leaves that a spanning
tree of $G$ can have among the internal vertices of $P_{uv}$.
Now we are ready to present the FPT algorithm in Algorithm~\ref{alg:FPT}.

\begin{algorithm}[htb]
\caption{An FPT algorithm for \maxleaf}\medskip
\label{alg:FPT}

INPUT: a \maxleaf\ instance $(G,k)$.\medskip

1) {\bf while} $G$ has a 2-terminal diamond or 2-terminal blossom subgraph {\bf do}

\q \q     $G:=$the result of applying (F1) or (F2) to $G$

\q \q     $k:=k-1$

\q {\bf end while}
   
2) {\bf if} $\nt(G)\geq 3k$ or $|L(G)|\geq k$ or $k\leq 2$ {\bf then} return(YES) {\bf endif}

3) construct $\H(G)$ and $c$

4) {\bf for} all $L\subseteq \Vt(G)$ with $|L|\leq k$ {\bf do}

\q \q     {\bf if} $G$ has a \ST\ $T$ with $L \subseteq L(T)$ and
$|L|+|L(T)\bs \Vt(G)|\geq k$ {\bf then} 

\q \q \q return(YES)

\q \q {\bf endif}

\q {\bf endfor}
       
5) return(NO)

\end{algorithm}

In the following proofs, we will use the fact that for any $L\subset
V(G)$, a \ST\ $T$ of $G$ with $L\subseteq L(T)$ exists if and only if
$G-L$ is connected and $V(G)\bs L$ is a dominating set.  
The decision in  Step 4 can be made in polynomial time in the size of $\H(G)$.
The essential step is to solve a minimum weight spanning tree
problem on $\H(G)-L$, using edge costs $c$. Lemma~\ref{lem:step4}
contains the details.

\newcounter{lem7}
\setcounter{lem7}{\value{lemma}}
\begin{lemma}\label{lem:step4}
\label{forloop_poly}
Let $(G,k)$ be a \maxleaf\ instance for which $\H(G)$ and $c$ are
non-empty and known. For any $L\subseteq \Vt(G)$, deciding whether $G$
has a \ST\ $T$ with $L\subseteq L(T)$ and $|L|+|L(T)\bs \Vt(G)|\geq k$
can be done in time polynomial in the size of $\H(G)$. 
\end{lemma}
\proof
Let $S=\H(G)$.
A \ST\ $T$ of $G$ with $L\subseteq L(T)$ exists if and only if
$V(G)\bs L$ is a connected, dominating set of $G$. This is the case if and
only if $V(S)\bs L$ is a connected, dominating set of $S$ and there is
no edge $uv\in E(S)$ with $u,v\in L$ and $c(uv)\geq 1$. These
properties can be checked in time polynomial in the size of $S$. 

Now suppose at least one \ST\ $T_0$ of $G$  with $L(T_0)\subseteq L$
exists. We show how to construct such a \ST\ $T$ that maximizes $|L(T)\bs \Vt(G)|$. 
This process is illustrated in Figure~\ref{fig:FPTtreeconstr}. White vertices indicate vertices in $L$.
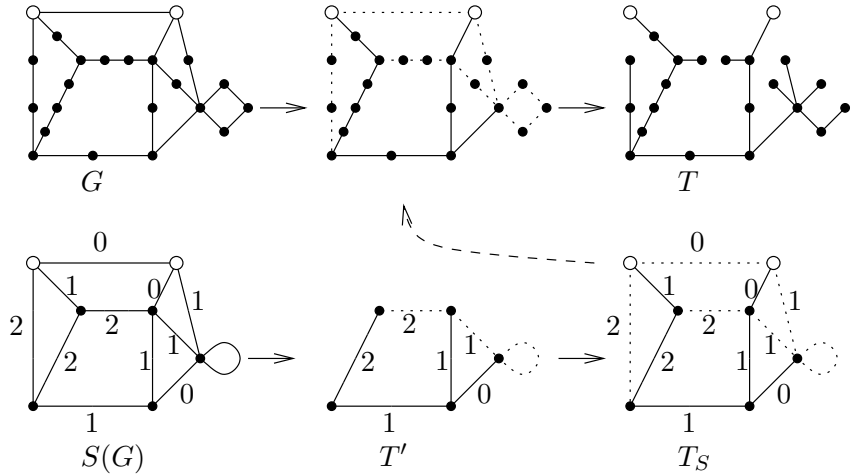
\PsFigCap{100}{FPTtreeconstr}{Constructing a tree $T$ with $L\subseteq L(T)$ using $S(G)$.}

First let $T'$ be a minimum weight \ST\ of $S-L$, with respect to the
cost function $c$. (This tree can be found in polynomial time.) 
A \ST\ $T_S$ of $S$ is obtained from $T'$ by connecting every  $u\in L$ to
$T'$ by an edge $e$ which has minimum cost $c(e)$. So $L\subseteq L(T_S)$.
 
A spanning tree $T$ of the subdivision $G$ of $S$ can be
obtained from $T_S$ in the following way. All edges from paths of $G$ that
correspond to edges in $T_S$ are in $T$. At this stage $T$ need not be
spanning. The inner vertices of a path $P_{uv}$ of $G$ corresponding to an edge
$uv\in E(S)\bs E(T_S)$ with  $u\in L, v\not\in L$ are connected to $T$ such that
$u$ remains a leaf. The inner vertices of a path $P_{uv}$ with 
$uv\in E(S)\bs E(T_S)$ and $u,v\not\in L$ can be connected to $T$ such that
$c(uv)$ of them become leaves of $T$. Finally, vertices of cycles of
$G$ that correspond to loops of $S$ are connected to $T$ such that two
of them become leaves. 
Note that an edge $uv\in E(S)\bs E(T_S)$ with $u,v\in L$ must correspond
to a path $P_{uv}$ in $G$ with no inner vertex since $G-L$ is connected. At this
point, every vertex of $G$ is connected to $T$, without introducing
cycles, hence $T$ is a \ST. 

It can be verified that $T$ was constructed such that $L(T)\bs \Vt(G)$
is maximized, under the condition that $L\subseteq L(T)$. For this it
is essential that $T'$ was chosen to be a minimal spanning tree of
$S-L$. We omit the formal proof of this fact here, noting that it is
similar to the proof given in~\cite{BBW03} and~\cite{Bon06}.  

Observe that $T$ does not actually have to be constructed for the
decision. Thus only $\H(G)$ needs to be considered, and the statement
follows.\qed 

\newcounter{lem8}
\setcounter{lem8}{\value{lemma}}
\begin{lemma}
\label{alg_correct}
Algorithm~\ref{alg:FPT} returns YES if and only if its input $(G,k)$
is a YES-instance for \maxleaf. 
\end{lemma}
\proof
Lemma~\ref{FPT_reds_correct} shows that it suffices to prove the
statement for the reduced instance $(G,k)$.
Note that the reduced instance $(G,k)$ is again simple, connected and
non-trivial and does not contain 2-terminal diamonds or 2-terminal
blossoms and thus also no \TEN s or \TEB s. So if $\nt(G)\geq 3k$,
then $(G,k)$ is a YES-instance by Theorem~\ref{thm:main}. The
correctness of the other cases in which the algorithm returns YES is
easily checked.  

Suppose that $(G,k)$ is a YES-instance. We show that the
algorithm indeed returns YES. If $\Vt(G)=\emptyset$, then
$G$ is a path or cycle, so $(G,k)$ being a 
YES-instance implies $k\leq 2$, and YES is returned in Step 2. 
Suppose $\Vt(G)\not=\emptyset$, Step
2 did not return YES, and $T$ is a \ST\ of $G$ with at least $k$
leaves. If $|L(T)\cap \Vt(G)|\geq k$, then some set $L\subseteq
L(T)\cap \Vt(G)$ with $|L|=k$ is 
considered in the algorithm. Clearly also a \ST\ $T'$ exists with
$L\subseteq L(T')$, so the algorithm returns YES in Step 4.  
On the other hand, if $L=L(T)\cap
\Vt(G)$ has fewer than $k$ elements, then $L$ itself will be
considered in Algorithm~\ref{alg:FPT}, and in addition $|L|+|L(T)\bs
\Vt(G)|=|L(T)|\geq k$. 
The algorithm will then return YES in Step 4.\qed 

Lemma~\ref{alg_correct} proves the correctness of Algorithm~1, the
claimed time complexity will be proved next. 
Together this proves Theorem~\ref{thm:FPT}. 
We repeat the statement for convenience.
\setcounter{repeatthm}{\value{thm2}}
\begin{repeatthm}
There exists an FPT algorithm for \maxleaf\ with time complexity
$O(m)+O^*(6.75^k)$, where $m$ denotes the size of the input graph and
$k$ the desired number of leaves.
\end{repeatthm}
\proof It only remains to prove the complexity bound. 

The first three steps
can be done in linear time by building the proper data structures. For
this it is essential that the degree of non-terminal vertices of
2-terminal diamonds and blossoms is bounded by a constant. We give
more details now.

Assume that $G$ is represented by doubly linked adjacency lists. 
It is possible to detect 2-terminal diamonds and 2-terminal blossoms
because all of their non-terminal vertices have degree at most~4.  
For every vertex $u$ of degree at most~4, we can store the position
of $u$ in the adjacency lists of its neighbors. To do this for every
such $u$, all adjacency lists of the graph need to be scanned only
once. This information now makes it possible to apply (F1) and (F2) in
constant time for each diamond and blossom respectively.  

Step 2 can obviously be done in linear time. For Step 3 we switch to a
representation using arrays for the edges, which contain the labels of
the two end vertices and the edge weights. We also store 
the vertex degrees, and the labels of
the incident edges for vertices of degree~2. This representation
allows us to do the following operations in constant time: suppressing
a degree~2 vertex, calculating the resulting edge weight, and
updating the representation. 

Thus Steps 1-3 can be performed in linear time and it only remains to
consider the complexity of Step~4.
Since the reductions in Step 1 do not increase the number of vertices or the
value of $k$, we may assume that $n$ and $k$
are the number of vertices and the parameter of the reduced instance, as
it is after Step 1. 

Step 4 of the algorithm is only executed when $|\Vt(G)|<3k$ and
$|L(G)|<k$. Furthermore $V(\H(G))=L(G)\cup \Vt(G)$, so every iteration
of the for-loop of Step 4 takes time polynomial in $k$
(Lemma~\ref{forloop_poly}). This for-loop is executed once for every
subset $L\subseteq \Vt(G)$ with $|L|\leq k$. Using $|\Vt(G)|\leq 3k$,
the number of such sets can be verified to be $O(k{3k \choose
  k})$. Using Stirling's approximation $x!\approx x^x e^{-x}
\sqrt{2\pi x}$, we obtain 
\[
{3k \choose k}=
\frac{(3k)!}{(2k)!k!}\in 
O\left(\frac{(3k)^{3k}}{e^{3k}}\cdot\frac{e^{2k}}{(2k)^{2k}}\cdot \frac{e^{k}}{k^k}\right)=
O\left(\frac{3^{3k}}{2^{2k}}\right)=O(6.75^k).
\]
This concludes the proof.\qed

We remark that we did not optimize the polynomial factor suppressed by the $O^*$ notation, but it can be seen to be a practical, low degree polynomial.

\section{Conclusions}
\label{sec:conclusions}

We conclude with some remarks about possible improvements. 
Theorem~\ref{thm:main}
can be strengthened at the cost of lengthier proofs. An extended version
of this paper will show that $\nt(G)/3+2$ leaves can be obtained
whenever $G$ is not cubic, and not equal to $\Gs$.
It can also be shown that in order to obtain this bound, \TEB s do not
have to be excluded; it suffices to only exclude larger structures
like the flower in Figure~\ref{fig:blossom+flowerX}. This way a
bound of $4\nt(G)/13+c$ can be proved when only \TEN s are excluded.

Besides optimizing the parameter function of FPT algorithms, 
another goal is to find better {\em 
kernelizations}. For \maxleaf, a kernelization can be defined as a
preprocessing method that reduces the input to an instance $(G,k)$
which is either a YES-instance, or has $|V(G)|\leq f(k)$ for some
function $f(k)$. The current best kernelization for \maxleaf\ has
$f(k)=3.75k$, see~\cite{EFL05}. Since our goal was to avoid
preprocessing as much as possible, our method does not give a
kernelization by itself. But for instance our approach can be combined
with three simple reduction rules from~\cite{BBW03} and~\cite{EFL05}
to remove all leaves and adjacent degree~2 vertices. This yields a
$7k$ kernelization: a reduced NO-instance has less than $3k$ vertices of degree at least 3 by Theorem~\ref{thm:main}, and less than $4k$ vertices of degree 2.

\bibliographystyle{plain}
\bibliography{leafytrees_medium}

\end{document}